\documentclass[12pt]{article}

\usepackage{amsmath,amssymb,amsthm,amscd,a4wide}
\usepackage[enableskew]{youngtab}

\begin{document}

\newcommand{\End}{{\rm{End}\ts}}
\newcommand{\Hom}{{\rm{Hom}}}
\newcommand{\Mat}{{\rm{Mat}}}
\newcommand{\ad}{{\rm{ad}\ts}}
\newcommand{\ch}{{\rm{ch}\ts}}
\newcommand{\chara}{{\rm{char}\ts}}
\newcommand{\ind}{{\rm{ind}\ts}}
\newcommand{\diag}{ {\rm diag}}
\newcommand{\pr}{^{\tss\prime}}
\newcommand{\non}{\nonumber}
\newcommand{\wt}{\widetilde}
\newcommand{\wh}{\widehat}
\newcommand{\ot}{\otimes}
\newcommand{\la}{\lambda}
\newcommand{\ls}{\ts\lambda\ts}
\newcommand{\La}{\Lambda}
\newcommand{\De}{\Delta}
\newcommand{\al}{\alpha}
\newcommand{\be}{\beta}
\newcommand{\ga}{\gamma}
\newcommand{\Ga}{\Gamma}
\newcommand{\ep}{\epsilon}
\newcommand{\ka}{\kappa}
\newcommand{\vk}{\varkappa}
\newcommand{\vt}{\vartheta}
\newcommand{\si}{\sigma}
\newcommand{\vp}{\varphi}
\newcommand{\de}{\delta}
\newcommand{\ze}{\zeta}
\newcommand{\om}{\omega}
\newcommand{\ee}{\epsilon^{}}
\newcommand{\su}{s^{}}
\newcommand{\hra}{\hookrightarrow}
\newcommand{\ve}{\varepsilon}
\newcommand{\ts}{\,}
\newcommand{\vac}{\mathbf{1}}
\newcommand{\di}{\partial}
\newcommand{\qin}{q^{-1}}
\newcommand{\tss}{\hspace{1pt}}
\newcommand{\Sr}{ {\rm S}}
\newcommand{\U}{ {\rm U}}
\newcommand{\BL}{ {\overline L}}
\newcommand{\BE}{ {\overline E}}
\newcommand{\BP}{ {\overline P}}
\newcommand{\AAb}{\mathbb{A}\tss}
\newcommand{\CC}{\mathbb{C}\tss}
\newcommand{\KK}{\mathbb{K}\tss}
\newcommand{\QQ}{\mathbb{Q}\tss}
\newcommand{\SSb}{\mathbb{S}\tss}
\newcommand{\ZZ}{\mathbb{Z}\tss}
\newcommand{\X}{ {\rm X}}
\newcommand{\Y}{ {\rm Y}}
\newcommand{\Z}{{\rm Z}}
\newcommand{\Ac}{\mathcal{A}}
\newcommand{\Lc}{\mathcal{L}}
\newcommand{\Mc}{\mathcal{M}}
\newcommand{\Pc}{\mathcal{P}}
\newcommand{\Qc}{\mathcal{Q}}
\newcommand{\Tc}{\mathcal{T}}
\newcommand{\Sc}{\mathcal{S}}
\newcommand{\Bc}{\mathcal{B}}
\newcommand{\Ec}{\mathcal{E}}
\newcommand{\Fc}{\mathcal{F}}
\newcommand{\Hc}{\mathcal{H}}
\newcommand{\Uc}{\mathcal{U}}
\newcommand{\Vc}{\mathcal{V}}
\newcommand{\Wc}{\mathcal{W}}
\newcommand{\Yc}{\mathcal{Y}}
\newcommand{\Ar}{{\rm A}}
\newcommand{\Br}{{\rm B}}
\newcommand{\Ir}{{\rm I}}
\newcommand{\Fr}{{\rm F}}
\newcommand{\Jr}{{\rm J}}
\newcommand{\Or}{{\rm O}}
\newcommand{\GL}{{\rm GL}}
\newcommand{\Spr}{{\rm Sp}}
\newcommand{\Rr}{{\rm R}}
\newcommand{\Zr}{{\rm Z}}
\newcommand{\gl}{\mathfrak{gl}}
\newcommand{\middd}{{\rm mid}}
\newcommand{\ev}{{\rm ev}}
\newcommand{\Pf}{{\rm Pf}}
\newcommand{\Norm}{{\rm Norm\tss}}
\newcommand{\oa}{\mathfrak{o}}
\newcommand{\spa}{\mathfrak{sp}}
\newcommand{\osp}{\mathfrak{osp}}
\newcommand{\g}{\mathfrak{g}}
\newcommand{\h}{\mathfrak h}
\newcommand{\n}{\mathfrak n}
\newcommand{\z}{\mathfrak{z}}
\newcommand{\Zgot}{\mathfrak{Z}}
\newcommand{\p}{\mathfrak{p}}
\newcommand{\sll}{\mathfrak{sl}}
\newcommand{\agot}{\mathfrak{a}}
\newcommand{\qdet}{ {\rm qdet}\ts}
\newcommand{\Ber}{ {\rm Ber}\ts}
\newcommand{\HC}{ {\mathcal HC}}
\newcommand{\cdet}{ {\rm cdet}}
\newcommand{\tr}{ {\rm tr}}
\newcommand{\gr}{ {\rm gr}}
\newcommand{\str}{ {\rm str}}
\newcommand{\loc}{{\rm loc}}
\newcommand{\Gr}{{\rm G}}
\newcommand{\sgn}{ {\rm sgn}\ts}
\newcommand{\ba}{\bar{a}}
\newcommand{\bb}{\bar{b}}
\newcommand{\bi}{\bar{\imath}}
\newcommand{\bj}{\bar{\jmath}}
\newcommand{\bk}{\bar{k}}
\newcommand{\bl}{\bar{l}}
\newcommand{\hb}{\mathbf{h}}
\newcommand{\Sym}{\mathfrak S}
\newcommand{\fand}{\quad\text{and}\quad}
\newcommand{\Fand}{\qquad\text{and}\qquad}
\newcommand{\For}{\qquad\text{or}\qquad}
\newcommand{\OR}{\qquad\text{or}\qquad}
\newcommand{\emp}{\mbox{}}
\newcommand{\atopn}[2]{\genfrac{}{}{0pt}{}{#1}{#2}}

\renewcommand{\theequation}{\arabic{section}.\arabic{equation}}

\newtheorem{thm}{Theorem}[section]
\newtheorem{lem}[thm]{Lemma}
\newtheorem{prop}[thm]{Proposition}
\newtheorem{cor}[thm]{Corollary}
\newtheorem{conj}[thm]{Conjecture}
\newtheorem*{thm-ref}{Theorem}
\newtheorem*{mthm}{Main Theorem}
\newtheorem*{mthma}{Theorem A}
\newtheorem*{mthmb}{Theorem B}

\theoremstyle{definition}
\newtheorem{defin}[thm]{Definition}

\theoremstyle{remark}
\newtheorem{remark}[thm]{Remark}
\newtheorem{example}[thm]{Example}

\newcommand{\bth}{\begin{thm}}
\renewcommand{\eth}{\end{thm}}
\newcommand{\bpr}{\begin{prop}}
\newcommand{\epr}{\end{prop}}
\newcommand{\ble}{\begin{lem}}
\newcommand{\ele}{\end{lem}}
\newcommand{\bco}{\begin{cor}}
\newcommand{\eco}{\end{cor}}
\newcommand{\bde}{\begin{defin}}
\newcommand{\ede}{\end{defin}}
\newcommand{\bex}{\begin{example}}
\newcommand{\eex}{\end{example}}
\newcommand{\bre}{\begin{remark}}
\newcommand{\ere}{\end{remark}}
\newcommand{\bcj}{\begin{conj}}
\newcommand{\ecj}{\end{conj}}

\newcommand{\bal}{\begin{aligned}}
\newcommand{\eal}{\end{aligned}}
\newcommand{\beq}{\begin{equation}}
\newcommand{\eeq}{\end{equation}}
\newcommand{\ben}{\begin{equation*}}
\newcommand{\een}{\end{equation*}}

\newcommand{\bpf}{\begin{proof}}
\newcommand{\epf}{\end{proof}}

\def\beql#1{\begin{equation}\label{#1}}

\title{\Large\bf Quantization of the shift of argument subalgebras\\ in type $A$}

\author{{Vyacheslav Futorny\quad and\quad Alexander Molev}}

\date{} 
\maketitle

\vspace{35 mm}

\begin{abstract}
Given a simple Lie algebra $\g$ and an element $\mu\in\g^*$,
the corresponding shift of argument subalgebra of $\Sr(\g)$ is
Poisson commutative.
In the case where $\mu$ is regular,
this subalgebra is known
to admit a quantization, that is, it can be lifted to a commutative subalgebra
of $\U(\g)$. We show that if $\g$ is of type $A$, then
this property extends to arbitrary $\mu$, thus proving a conjecture
of Feigin, Frenkel and Toledano Laredo. The proof relies on an explicit
construction of generators of the center of the affine vertex algebra
at the critical level.

%

\end{abstract}


\vspace{45 mm}

\noindent
Department of Mathematics,
University of S\~{a}o Paulo,\newline
Caixa Postal 66281, S\~{a}o Paulo, SP 05315-970, Brazil\newline
vfutorny@gmail.com

\vspace{7 mm}

\noindent
School of Mathematics and Statistics\newline
University of Sydney,
NSW 2006, Australia\newline
alexander.molev@sydney.edu.au

%

\newpage

\section{Introduction}
\label{sec:int}
\setcounter{equation}{0}

\paragraph{Shift of argument subalgebras.}

Let $\g$ be a simple Lie algebra over $\CC$ with basis
elements $Y_1,\dots,Y_l$ and the corresponding structure constants
$c_{ij}^{\tss k}$.
The symmetric algebra $\Sr(\g)$ can be equipped with
the {\it Lie--Poisson bracket\/} defined on the elements
of the Lie algebra by
\beql{liepoisson}
\{Y_i,Y_j\}=\sum_{k=1}^l\tss c_{ij}^{\tss k}\ts Y_k.
\eeq

Let $P=P(Y_1,\dots,Y_l)$ be an element of $\Sr(\g)$ of a certain degree $d$.
Fix any element $\mu\in\g^*$ and let $z$ be a variable.
Make the substitution $Y_i\mapsto Y_i+z\ts\mu(Y_i)$
and expand as a polynomial in $z$,
\ben
P\big(Y_1+z\ts\mu(Y_1),\dots,Y_l+z\ts\mu(Y_l)\big)
=P^{(0)}+P^{(1)} z+\dots+P^{(d)} z^d
\een
to define elements $P^{(i)}\in \Sr(\g)$ associated with $P$ and $\mu$.
Denote by
$\overline\Ac_{\mu}$ the subalgebra
of $\Sr(\g)$ generated by all elements $P^{(i)}$
associated with all $\g$-invariants $P\in \Sr(\g)^{\g}$.
The subalgebra $\overline\Ac_{\mu}$ of $\Sr(\g)$ is known
as the {\it Mishchenko--Fomenko subalgebra\/} or
{\it shift of argument subalgebra\/}. Its key property observed in \cite{mf:ee} states
that $\overline\Ac_{\mu}$ is Poisson commutative; that is, $\{R,S\}=0$
for any elements $R,S\in\overline\Ac_{\mu}$.

We will identify $\g^*$ with $\g$ via a symmetric invariant bilinear form
(see \eqref{killi} below) and let $n$ denote the rank of $\g$.
An element $\mu\in\g^*\cong\g$ is called {\em regular},
if the centralizer $\g^{\mu}$ of $\mu$ in $\g$ has minimal possible
dimension; this minimal dimension coincides with $n$.
The subalgebra $\Sr(\g)^{\g}$ admits
a family $P_1,\dots,P_n$ of algebraically independent generators
of respective degrees $d_1,\dots,d_n$. If the element $\mu\in\g^*$ is regular,
then $\overline\Ac_{\mu}$ has the properties:
\begin{enumerate}
\item[$i)$]
the subalgebra $\overline\Ac_{\mu}$
of $\Sr(\g)$ is maximal Poisson commutative;
\item[$ii)$]\label{pki}
the elements
$P^{(i)}_{k}$ with $k=1,\dots,n$ and $i=0,1,\dots,d_k-1$,
are algebraically independent generators of $\overline\Ac_{\mu}$.
\end{enumerate}

Property $i)$ is a theorem of Panyushev and Yakimova~\cite{py:as}; the case of
regular semisimple $\mu$ is due to Tarasov~\cite{t:ms}.
Property $ii)$ is due to Bolsinov~\cite{b:cf}; the regular semisimple case
goes back to the original paper \cite{mf:ee}. Another proof of $ii)$
was given in~\cite{fft:gm}.

\paragraph{Vinberg's problem.}

The universal enveloping algebra $\U(\g)$ is equipped with a canonical filtration and the
associated graded algebra $\gr\ts\U(\g)$ is isomorphic to
$\Sr(\g)$.
Given that the subalgebra
$\overline\Ac_{\mu}$ of $\Sr(\g)$ is Poisson commutative,
one could look for
a commutative
subalgebra $\Ac_{\mu}$ of $\U(\g)$ which ``quantizes" $\overline\Ac_{\mu}$
in the sense that $\gr\ts\Ac_{\mu}=\overline\Ac_{\mu}$.
This quantization problem was raised by Vinberg in \cite{v:sc},
where, in particular, some commuting families of elements of $\U(\g)$
were produced. A positive solution of Vinberg's problem
was given by Rybnikov~\cite{r:si} (for regular semisimple $\mu$)
and Feigin, Frenkel and Toledano Laredo~\cite{fft:gm} (for any regular $\mu$)
with the use of the center of the associated affine vertex algebra at the critical level
(also known as the {\em Feigin--Frenkel center}).
To briefly outline
the solution, equip $\g$ with
a standard symmetric invariant bilinear form $\langle\ts\ts,\ts\rangle$
defined as the normalized Killing form
\beql{killi}
\langle X,Y\rangle=\frac{1}{2\tss h^{\vee}}\ts\tr\ts\big(\ad\tss X\ts\ad\tss Y\big),
\eeq
where $h^{\vee}$ is the {\it dual Coxeter number\/} for $\g$.
The corresponding {\it affine Kac--Moody algebra\/} $\wh\g$
is the central
extension
\beql{km}
\wh\g=\g\tss[t,t^{-1}]\oplus\CC K,
\eeq
where $\g[t,t^{-1}]$ is the Lie algebra of Laurent
polynomials in $t$ with coefficients in $\g$. For any $r\in\ZZ$ and $X\in\g$
we set $X[r]=X\ts t^r$. The commutation relations of the Lie algebra $\wh\g$
have the form
\ben
\big[X[r],Y[s]\big]=[X,Y][r+s]+r\ts\de_{r,-s}\langle X,Y\rangle\ts K,
\qquad X, Y\in\g,
\een
and the element $K$ is central in $\wh\g$.
For any $\ka\in\CC$ denote by $\U_{\ka}(\wh\g)$ the quotient of $\U(\wh\g)$
by the ideal generated by $K-\ka$. The value $\ka=-h^{\vee}$
corresponds to the {\it critical level\/}.
Let $\Ir$ denote the left ideal of $\U_{-h^{\vee}}(\wh\g)$ generated by $\g[t]$
and let $\Norm\tss\Ir$ be its normalizer,
\ben
\Norm\tss\Ir=\{v\in \U_{-h^{\vee}}(\wh\g)\ |\ \Ir\tss v\subseteq \Ir\}.
\een
The normalizer is a subalgebra of $\U_{-h^{\vee}}(\wh\g)$, and $\Ir$
is a two-sided ideal of $\Norm\tss\Ir$.
The {\it Feigin--Frenkel center\/} $\z(\wh\g)$ is the associative algebra
defined as the quotient
\beql{ffnorm}
\z(\wh\g)=\Norm\tss\Ir/\Ir.
\eeq
By the Poincar\'e--Birkhoff--Witt theorem, the quotient of the algebra
$\U_{-h^{\vee}}(\wh\g)$ by the left ideal $\Ir$
is isomorphic to the universal enveloping algebra
$\U\big(t^{-1}\g[t^{-1}]\big)$, as a vector space. Hence, we have
a vector space embedding
\ben
\z(\wh\g)\hra \U\big(t^{-1}\g[t^{-1}]\big).
\een
Since $\U\big(t^{-1}\g[t^{-1}]\big)$ is a subalgebra of $\U_{-h^{\vee}}(\wh\g)$,
the embedding
is an algebra homomorphism so that the Feigin--Frenkel center
$\z(\wh\g)$ can be regarded as a
subalgebra of $\U\big(t^{-1}\g[t^{-1}]\big)$. In fact, this subalgebra is {\it commutative\/}
which is not immediate from the definition, but can be seen by
identifying $\z(\wh\g)$ with the center of the affine vertex algebra
at the critical level. Furthermore,
by a theorem of Feigin and Frenkel~\cite{ff:ak} (see \cite{f:lc}
for a detailed exposition),
there exist elements $S_1,\dots,S_n\in \z(\wh\g)$ such that
\beql{genz}
\z(\wh\g)=\CC[T^{\tss r}S_l\ |\ l=1,\dots,n,\ \ r\geqslant 0],
\eeq
where $T$ is the derivation of the algebra $\U\big(t^{-1}\g[t^{-1}]\big)$
which is determined by the property that its commutator
with the operator of left multiplication by $X[r]$ is found by
\ben
\big[T,X[r]\big]=-r\tss X[r - 1],\qquad X\in\g,\quad r<0.
\een
We will call such family $S_1,\dots,S_n$ a {\em complete set of
Segal--Sugawara vectors for $\g$}. Another derivation $D$
of the algebra $\U\big(t^{-1}\g[t^{-1}]\big)$
is determined by the property
\ben
\big[D,X[r]\big]=-r\tss X[r],\qquad X\in\g,\quad r<0;
\een
and $D$ defines a grading on $\U\big(t^{-1}\g[t^{-1}]\big)$.

Given any element $\mu\in\g^*$ and
a nonzero $z\in\CC$, the mapping
\beql{evalr}
\varrho^{}_{\ts\mu,z}:\U\big(t^{-1}\g[t^{-1}]\big)\to \U(\g),
\qquad X[r]\mapsto X z^r+\de_{r,-1}\ts\mu(X),\quad X\in\g,
\eeq
defines an algebra homomorphism. The image of $\z(\wh\g)$
under $\varrho^{}_{\ts\mu,z}$
is a commutative subalgebra of $\U(\g)$. It does not depend
on $z$ and is denoted by $\Ac_{\mu}$.
If $S\in \U\big(t^{-1}\g[t^{-1}]\big)$ is an element
of degree $d$ with respect to the grading defined by $D$,
then regarding $\varrho^{}_{\ts\mu,z}(S)$
as a polynomial in $z^{-1}$, define the elements $S^{(i)}\in\U(\g)$
by the expansion
\beql{rhoexp}
\varrho^{}_{\ts\mu,z}(S)=S^{(0)} z^{-d}+\dots+S^{(d-1)} z^{-1}+S^{(d)}.
\eeq
If $\mu\in\g^*$ is regular then the following holds:
\begin{enumerate}
\item[$i)$]
the subalgebra $\Ac_{\mu}$ of
$\U(\g)$ is maximal commutative;
\item[$ii)$]
if $S_1,\dots,S_n\in \z(\wh\g)$ are elements of
the respective degrees $d_1,\dots,d_n$ satisfying \eqref{genz},
then the elements
$S^{(i)}_{k}$ with $k=1,\dots,n$ and $i=0,1,\dots,d_k-1$
are algebraically independent generators of $\Ac_{\mu}$;
\item[$iii)$]
$\gr\ts\Ac_{\mu}=\overline\Ac_{\mu}$.
\end{enumerate}
This is derived with the use of the respective properties of the
algebra $\overline\Ac_{\mu}$; see \cite{fft:gm} for proofs.
The subalgebra $\Ac_{\mu}$ was further studied in \cite{ffr:oi}
where its spectra in finite-dimensional
irreducible representations of $\g$ were described.

Note that both algebras $\Ac_{\mu}$ and $\overline\Ac_{\mu}$ are defined
for arbitrary elements $\mu\in\g^*$. Given that the property $iii)$ holds
for all regular $\mu$, it was conjectured in \cite[Conjecture~1]{fft:gm},
that this property is valid for all $\mu$. As a consequence of our
main result, we obtain a proof of
this conjecture for type $A$; see the Main Theorem below.
In particular, this gives another proof of $iii)$ for regular $\mu$.
More precisely, we will work
with the reductive Lie algebra $\g=\gl_n$ and consider
the respective subalgebras $\overline\Ac_{\mu}\subset \Sr(\gl_n)$
and $\Ac_{\mu}\subset \U(\gl_n)$.
The proof will be based on the use of explicit formulas for generators of
$\Ac_{\mu}$.

\paragraph{Generators of $\Ac_{\mu}$.}
For the Lie algebras $\g$ of type $A$, a few families of
explicit generators $S_1,\dots,S_n$ of $\z(\wh\g)$,
and hence generators of the subalgebra $\Ac_{\mu}$, were
produced by Chervov and Talalaev~\cite{ct:qs} by extending Talalaev's
work \cite{t:qg}; see also \cite{cm:ho} and
\cite{mr:mm} where more direct proofs were given.
In types $B$, $C$ and $D$ such explicit generators were constructed in \cite{m:ff}.
Note also earlier work of Nazarov and Olshanski~\cite{no:bs}, where
maximal commutative subalgebras of $\U(\g)$ were produced with the use
of Yangians; they quantize the Poisson algebras $\overline\Ac_{\mu}$
in all classical types for the case of regular semisimple $\mu$.
In a different form, a quantization of $\overline\Ac_{\mu}$
in type $A$ was provided by Tarasov~\cite{t:cs} via a symmetrization map.

We will work with a particular family of generators
of $\z(\wh\gl_n)$
which we recall below in Sec.~\ref{sec:gff}.
They allow us to define the associated family of generators
$\phi^{(k)}_m$ with $m=1,\dots,n$ and $k=0,\dots,m-1$
of the subalgebra $\Ac_{\mu}\subset\U(\gl_n)$; see \eqref{defvp} below.
There generators are algebraically independent if $\mu$ is
regular.

Our main result provides a way to choose
an algebraically independent family of generators $\phi^{(k)}_m$ of $\Ac_{\mu}$
for an arbitrary element $\mu$.
To describe
this subset, we will identify $\gl^*_n$ with $\gl_n$
via a symmetric bilinear form and regard $\mu$ as an $n\times n$ matrix.
Suppose that the distinct eigenvalues of $\mu$ are $\la_1,\dots,\la_r$
and the Jordan canonical form of $\mu$ is the direct sum of the respective Jordan blocks
$J_{\al^{(i)}_j}(\la_i)$ of sizes
$\al^{(i)}_1\geqslant \al^{(i)}_2\geqslant\dots  \geqslant \al^{(i)}_{s_i}\geqslant 1$.
We let $\al^{(i)}$ denote the corresponding Young diagram whose $j$-th row
is $\al^{(i)}_j$ and let $|\al^{(i)}|$ be the number of boxes of $\al^{(i)}$.
Given these data, introduce another Young diagram
$\ga=(\ga_1,\ga_2,\dots)$ by setting
\beql{gal}
\ga_l=\sum_{i=1}^r\sum_{j\geqslant l+1} \al^{(i)}_j,
\eeq
so that $\ga_l$ is the total number of boxes
which are strictly below the $l$-th rows in all diagrams $\al^{(i)}$.
Furthermore, associate the elements of the family $\phi^{(k)}_m$
with boxes of the diagram $\Ga=(n,n-1,\dots,1)$ so that
the $(i,j)$ box of $\Ga$ corresponds to $\phi^{(n-i-j+1)}_{n-j+1}$,
as illustrated:
\beql{Ga}
\Ga\quad=\qquad
\begin{matrix}
\phi^{(n-1)}_n & \phi^{(n-2)}_{n-1} & \dots & \phi^{(1)}_2 & \phi^{(0)}_1\\[0.5em]
\phi^{(n-2)}_n & \phi^{(n-3)}_{n-1} & \dots & \phi^{(0)}_2 & \\
\dots & \dots & \dots & & \\[0.5em]
\phi^{(1)}_n\phantom{-} & \phi^{(0)}_{n-1}\phantom{-} & & &\\[0.5em]
\phi^{(0)}_n\phantom{-} & & & &
\end{matrix}
\eeq
Note that the diagram $\ga$ is contained in $\Ga$.
We can now state our main theorem, where $\mu$ is an arbitrary
element of $\gl_n$ and $\Ga/\ga$ is the associated skew diagram.

\begin{mthm}
The elements $\phi^{(k)}_m$ corresponding to the boxes
of the skew diagram $\Ga/\ga$ are algebraically independent generators
of the subalgebra $\Ac_{\mu}$.
Moreover,
the subalgebra $\Ac_{\mu}$ is a quantization of $\overline\Ac_{\mu}$
so that $\gr\ts\Ac_{\mu}=\overline\Ac_{\mu}$.
\end{mthm}

By considering some other complete sets of Segal--Sugawara vectors,
we also show that the first part of the Main Theorem remains valid
if the elements $\phi^{(k)}_m$ are replaced with those of
other families; see Corollaries~\ref{cor:othgen} and \ref{cor:othgentwo} below.

\bex\label{ex:gensk}
Take $n=6$ and let $\mu$ be a nilpotent matrix with the Jordan blocks
of sizes $(2,2,1,1)$. Then
$\ga=(4,2,1)$ and the skew diagram $\Ga/\ga$ is
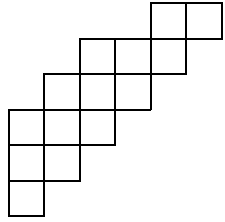
\ben
\young(::::\emp\emp,::\emp\emp\emp,:\emp\emp\emp,\emp\emp\emp,\emp\emp,\emp)
\een
so that the algebraically independent
generators of $\Ac_{\mu}$ are those corresponding to the boxes of $\Ga$,
excluding $\phi^{(2)}_3$, $\phi^{(3)}_4$, $\phi^{(3)}_5$, $\phi^{(4)}_5$,
$\phi^{(3)}_6$, $\phi^{(4)}_6$ and $\phi^{(5)}_6$.
\qed
\eex

Note also two extreme cases.
If $\mu$ is regular, then all Jordan blocks correspond to distinct
eigenvalues so that each $\al^{(i)}$ is a singe row diagram. Therefore,
$\ga=\varnothing$, so that all generators $\phi^{(k)}_m$
associated with the boxes of $\Ga$ are algebraically independent.
On the other hand,
for scalar matrices
$\mu$ we have $\ga=(n-1,n-2,\dots,1)$. In this case, $\Ac_{\mu}$ is generated
by $\phi^{(0)}_1,\dots,\phi^{(0)}_n$ and it coincides with the center of $\U(\gl_n)$.

Our proofs rely on {\em Bolsinov's completeness criterion} \cite[Theorem~3.2]{b:cf}
which applies to the shift of argument subalgebras associated
with an arbitrary Lie algebra $\g$. The required condition
for reductive Lie algebras is the equality
\beql{beconj}
\ind\g=\ind\g^{\mu}
\eeq
of the indices of $\g$ and the centralizer $\g^{\mu}$ of $\mu$ in $\g$, where
the {\em index} of an arbitrary Lie algebra
$\g$ is the minimal dimension of the stabilizers
$\g^x$, $x\in\g^*$, for the coadjoint representation.
In the case $\g=\gl_n$ and arbitrary $\mu\in\g$
this equality was claimed to be verified by Bolsinov
\cite[Sec.~3]{b:cp} (and was suggested to be extendable to
arbitrary semisimple Lie algebras)
and by Elashvili (private communication), but details were not published.
The first published proof is due to Yakimova~\cite{y:ic}, which extends to
all classical Lie algebras. The equality \eqref{beconj} is widely referred to
as the {\em Elashvili conjecture}, but should rather be called
the {\em Bolsinov--Elashvili conjecture}\footnote{
A.~Elashvili kindly informed us that the conjectural equality had
emerged from A.~Bolsinov's questions to him and so
it should also be attributed to the author
of \cite{b:cf}.}; see
e.g. \cite{cm:ic} for its proof covering all simple Lie algebras
and more references.

\smallskip

We are grateful to Alexey Bolsinov,
Alexander Elashvili,
Leonid Rybnikov and Alexander Veselov
for useful discussions.
The first author  was supported in part by the
CNPq grant (301320/2013-6) and by the
Fapesp grant (2010/50347-9).
This work was completed during the second author's visit to
the University of S\~{a}o Paulo. He would like to thank the
Department of Mathematics for the warm hospitality.

\section{Generators of $\z(\wh\gl_n)$}
\label{sec:gff}
\setcounter{equation}{0}

For $i,j\in\{1,\dots,n\}$ we will denote by $E_{ij}$
the standard basis
elements of $\gl_n$.
We extend the form
\eqref{killi} to the
invariant symmetric bilinear form on $\gl_n$ which is given by
\ben
\langle X,\ts Y\rangle=\tr\tss(X\tss Y)-\frac{1}{n}\ts\tr\tss X
\ts\tr\tss Y,\qquad X, Y\in\gl_n,
\een
where $X$ and $Y$ are regarded as $n\times n$ matrices.
Note that the kernel of the form is spanned by
the element $E_{11}+\dots+E_{nn}$, and its restriction to the subalgebra
$\sll_n$ is given by
\ben
\langle X,\ts Y\rangle=\tr\tss(X\tss Y),\qquad X, Y\in\sll_n.
\een
The affine Kac--Moody algebra $\wh\gl_n=\gl_n[t,t^{-1}]\oplus\CC K$
has the commutation relations
\beql{commrel}
\big[E_{ij}[r],E_{kl}[s\tss]\tss\big]
=\de_{kj}\ts E_{i\tss l}[r+s\tss]
-\de_{i\tss l}\ts E_{kj}[r+s\tss]
+r\tss\de_{r,-s}\ts K\Big(\de_{kj}\tss\de_{i\tss l}
-\frac{\de_{ij}\tss\de_{kl}}{n}\Big),
\eeq
and the element $K$ is central. The
critical level $-n$ coincides with
the negative of the dual Coxeter number for $\sll_n$.
We will work with
the extended Lie algebra $\wh\gl_n\oplus\CC\tau$ where the
additional element $\tau$ satisfies
the commutation relations
\beql{taur}
\big[\tau,X[r]\tss\big]=-r\ts X[r-1],\qquad
\big[\tau,K\big]=0.
\eeq
For any $r\in\ZZ$ combine the elements $E_{ij}[r]$
into the matrix $E[r]$ so that
\beql{matrer}
E[r]=\sum_{i,j=1}^n e_{ij}\ot E_{ij}[r]\in \End\CC^n\ot \U,
\eeq
where the $e_{ij}$ are the standard matrix units and
$\U$ stands
for the universal enveloping algebra of
$\wh\gl_n\oplus\CC\tau$.
For each $a\in\{1,\dots,m\}$
introduce the element $E[r]_a$ of the algebra
\beql{tenprka}
\underbrace{\End\CC^{n}\ot\dots\ot\End\CC^{n}}_m{}\ot\U
\eeq
by
\beql{matnota}
E[r]_a=\sum_{i,j=1}^{n}
1^{\ot(a-1)}\ot e_{ij}\ot 1^{\ot(m-a)}\ot E_{ij}[r].
\eeq
We let $H^{(m)}$ and $A^{(m)}$ denote the respective images of the
symmetrizer $h^{(m)}$ and anti-symmetrizer $a^{(m)}$ in the group algebra
for the symmetric group $\Sym_m$ under
its natural action on $(\CC^{n})^{\ot m}$. The elements $h^{(m)}$ and $a^{(m)}$
are the idempotents in the group algebra $\CC[\Sym_m]$ defined by
\ben
h^{(m)}=\frac{1}{m!}\ts\sum_{s\in\Sym_m} s
\Fand a^{(m)}=\frac{1}{m!}\ts\sum_{s\in\Sym_m} \sgn s\cdot s.
\een
We will identify $H^{(m)}$ and $A^{(m)}$ with the respective elements
$H^{(m)}\ot 1$ and $A^{(m)}\ot 1$ of the algebra \eqref{tenprka}.
Define the elements
$\phi^{}_{m\tss a},\psi^{}_{m\tss a}, \theta^{}_{m\tss a}\in
\U\big(t^{-1}\gl_n[t^{-1}]\big)$
by the expansions
\begin{align}\label{deftra}
\tr_{1,\dots,m}\ts A^{(m)} \big(\tau+E[-1]_1\big)\dots \big(\tau+E[-1]_m\big)
&=\phi^{}_{m\tss0}\ts\tau^m+\phi^{}_{m\tss1}\ts\tau^{m-1}
+\dots+\phi^{}_{m\tss m},\\[0.7em]
\label{deftrh}
\tr_{1,\dots,m}\ts H^{(m)} \big(\tau+E[-1]_1\big)\dots \big(\tau+E[-1]_m\big)
&=\psi^{}_{m\tss0}\ts\tau^m+\psi^{}_{m\tss1}\ts\tau^{m-1}
+\dots+\psi^{}_{m\tss m},
\end{align}
where
the traces are taken with respect to all $m$ copies of $\End\CC^n$
in \eqref{tenprka},
and
\beql{deftracepa}
\tr\ts \big(\tau+E[-1]\big)^m=\theta^{}_{m\tss0}\ts\tau^m+\theta^{}_{m\tss1}\ts\tau^{m-1}
+\dots+\theta^{}_{m\tss m}.
\eeq
Expressions like $\tau+E[-1]$ are understood as matrices, where
$\tau$ is regarded as the scalar matrix $\tau\tss 1$.
Furthermore, introduce the {\em column-determinant} of the matrix $\tau+E[-1]$
by
\beql{cdet}
\cdet\ts \big(\tau+E[-1]\big)=\sum_{\si\in\Sym_n} \sgn\si\cdot
\big(\tau+E[-1]\big)_{\si(1)\tss 1}\dots
\big(\tau+E[-1]\big)_{\si(n)\tss n}
\eeq
and expand it as a polynomial in $\tau$,
\beql{coldetal}
\cdet\ts \big(\tau+E[-1]\big)=\tau^n+\phi^{}_{1}\ts\tau^{n-1}
+\dots+\phi^{}_{n},\qquad \phi^{}_{m}\in
\U\big(t^{-1}\gl_n[t^{-1}]\big).
\eeq
We have the expansion of the noncommutative characteristic polynomial,
\beql{idtrd}
\cdet\ts \big(u+\tau+E[-1]\big)=\sum_{m=0}^n u^{n-m}\ts
\tr_{1,\dots,m}\ts A^{(m)} \big(\tau+E[-1]_1\big)\dots \big(\tau+E[-1]_m\big),
\eeq
where $u$ is a variable.
This implies the relations
\beql{chaam}
\phi^{}_{m\tss a}=\binom{n-a}{m-a}\ts\phi^{}_{a},\qquad
0\leqslant a\leqslant m\leqslant n.
\eeq
In particular, $\phi^{}_{m\tss m}=\phi^{}_{m}$ for $m=1,\dots,n$.

\bth\label{thm:allff}
All elements $\phi^{}_{m}$,
$\psi^{}_{m\tss a}$ and $\theta^{}_{m\tss a}$
belong to the Feigin--Frenkel center $\z(\wh\gl_n)$.
Moreover, each of the families
\ben
\phi^{}_{1},\dots,\phi^{}_{n},\qquad \psi^{}_{1\tss 1},\dots,\psi^{}_{n\tss n}
\Fand \theta^{}_{1\tss 1},\dots,\theta^{}_{n\tss n}
\een
is a complete set of Segal--Sugawara vectors for $\gl_n$.
\qed
\eth

This theorem goes back to \cite{ct:qs}, where the elements $\phi^{}_{m}$
were first discovered (in a slightly different form).
A direct proof of the theorem was given in \cite{cm:ho}.
The elements $\psi^{}_{m\tss a}$ are related to $\phi^{}_{m\tss a}$
through the quantum MacMahon Master Theorem of \cite{glz:qm}, while
a relationship between the $\phi^{}_{m\tss a}$ and $\theta^{}_{m\tss a}$
is provided by a Newton-type identity given in \cite[Theorem~15]{cfr:ap}.
Note that super-versions of these relations between the families
of Segal--Sugawara vectors for the Lie
superalgebra $\gl_{m|n}$ were given in the paper \cite{mr:mm}, which also provides
simpler arguments in the purely even case.

\section{Generators of $\Ac_{\mu}$}
\label{sec:gamu}
\setcounter{equation}{0}

In accordance with the results which we recalled in the Introduction,
the application of the homomorphism \eqref{evalr} to elements of $\z(\wh\gl_n)$
provided by Theorem~\ref{thm:allff} yields the corresponding families
of elements of the subalgebra $\Ac_{\mu}\subset\U(\gl_n)$ through the expansion
\eqref{rhoexp}. To give explicit formulas,
we will use the tensor
product algebra \eqref{tenprka}, where $\U$ will now denote the
algebra of differential operators whose elements are
finite sums of the form
\ben
\sum_{k,l\geqslant 0} u^{}_{kl}\ts z^{-k}\tss \di^{\ts l}_z,
\qquad u^{}_{kl}\in\U(\gl_n).
\een
Note that $\di_z$ emerges here as the image of the element $-\tau$ under
the extension of the homomorphism \eqref{evalr}.
As in \eqref{matrer}, we set
\ben
E=\sum_{i,j=1}^n e_{ij}\ot E_{ij}\in \End\CC^n\ot \U(\gl_n),
\een
and extend the notation \eqref{matnota} to the matrices $E$, $\mu$
and $M=-\di_z+\mu+Ez^{-1}$.
Assuming that $\mu\in\gl_n$ is arbitrary,
introduce the polynomials $\phi^{}_{m\tss a}(z)$, $\psi^{}_{m\tss a}(z)$
and $\theta^{}_{m\tss a}(z)$
in $z^{-1}$ (depending on $\mu$) with coefficients in $\U(\gl_n)$
by the expansions
\ben
\bal
\tr_{1,\dots,m}\ts A^{(m)} M_1\dots
M_m&=\phi^{}_{m\tss0}(z)\ts\di_z^{\ts m}+\phi^{}_{m\tss1}(z)\ts\di_z^{\ts m-1}
+\dots+\phi^{}_{m\tss m}(z),\\[0.5em]
\tr_{1,\dots,m}\ts H^{(m)} M_1\dots
M_m&=\psi^{}_{m\tss0}(z)\ts\di_z^{\ts m}+\psi^{}_{m\tss1}(z)\ts\di_z^{\ts m-1}
+\dots+\psi^{}_{m\tss m}(z),
\eal
\een
and
\ben
\tr\ts M^m
=\theta^{}_{m\tss0}(z)\ts\di_z^{\ts m}+\theta^{}_{m\tss1}(z)\ts\di_z^{\ts m-1}
+\dots+\theta^{}_{m\tss m}(z).
\een
Furthermore, following \eqref{coldetal}
define the polynomials $\phi^{}_{a}(z)$ by expanding
the column-determinant
\beql{cdetcoma}
\cdet\ts M=\phi^{}_{0}(z)\ts\di_z^{\ts n}+\phi^{}_{1}(z)
\ts\di_z^{\ts n-1}
+\dots+\phi^{}_{n}(z).
\eeq
By \eqref{chaam} we have
\ben
\phi^{}_{m\tss a}(z)=\binom{n-a}{m-a}\ts\phi^{}_{a}(z),\qquad
0\leqslant a\leqslant m\leqslant n,
\een
and so
$\phi^{}_{m\tss m}(z)=\phi^{}_{m}(z)$ for all $m$.
Introduce the coefficients of polynomials by
\ben
\bal
\phi^{}_{m}(z)&=\phi^{\tss(0)}_{m}z^{-m}
+\dots+\phi^{\tss(m-1)}_{m}z^{-1}
+\phi^{\tss(m)}_{m},\\[0.5em]
\psi^{}_{m\tss m}(z)&=\psi^{\tss(0)}_{m\tss m}z^{-m}
+\dots+\psi^{\tss(m-1)}_{m\tss m}z^{-1}
+\psi^{\tss(m)}_{m\tss m},
\eal
\een
and
\ben
\theta^{}_{m\tss m}(z)=\theta^{\tss(0)}_{m\tss m}z^{-m}
+\dots+\theta^{\tss(m-1)}_{m\tss m}z^{-1}
+\theta^{\tss(m)}_{m\tss m}.
\een

By Theorem~\ref{thm:allff} and the general results
of \cite{fft:gm} and \cite{r:si} we get the following.

\bth\label{thm:comsera}
Given any $\mu\in\gl_n$, all coefficients of the polynomials
$\phi^{}_{m}(z)$, $\psi^{}_{m\tss a}(z)$ and
$\theta^{}_{m\tss a}(z)$
belong to the commutative subalgebra $\Ac_{\mu}$ of $\U(\gl_n)$.
Moreover, the elements of each
of the families
\ben
\phi^{\tss(k)}_{m},\qquad \psi^{\tss(k)}_{m\tss m}\Fand
\theta^{\tss(k)}_{m\tss m}
\een
with $m=1,\dots,n$ and $k=0,1,\dots,m-1$, are generators of the algebra $\Ac_{\mu}$.
If $\mu$ is regular, then each of these families is algebraically
independent.
\qed
\eth

\bex\label{ex:tracom}
Using the family $\theta^{\tss(k)}_{m\tss m}$ we get the following
algebraically
independent generators of the algebra $\Ac_{\mu}$ for regular $\mu$:
\ben
\bal
\text{for}\quad\gl_2: &\qquad\tr\ts E,\quad \tr\ts \mu\tss E,\quad\tr\ts E^2\\[0.3em]
\text{for}\quad\gl_3: &\qquad\tr\ts E,\quad \tr\ts \mu\tss E,\quad\tr\ts \mu^2 E,
\quad \tr\ts E^2,\quad \tr\ts \mu\tss E^2,\quad \tr\ts E^3\\[0.3em]
\text{for}\quad\gl_4: &\qquad\tr\ts E,\quad \tr\ts \mu\tss E,\quad\tr\ts \mu^2 E,
\quad\tr\ts \mu^3 E,
\quad \tr\ts E^2,\quad \tr\ts \mu\tss E^2,\\[0.3em]
&\qquad\qquad 2\ts\tr\ts \mu^2 E^2+\tr\ts(\mu\tss E)^2,\quad
\tr\ts E^3,\quad\tr\ts \mu\tss E^3,\quad \tr\ts E^4.
\eal
\een
\eex

\section{Proof of the Main Theorem}
\label{sec:ptab}
\setcounter{equation}{0}

Note that $M=-\di_z+\mu+Ez^{-1}$ is a {\em Manin matrix} and therefore the polynomials
$\phi^{}_{m\tss a}(z)$ and $\psi^{}_{m\tss a}(z)$ admit expressions in terms of
noncommutative minors and permanents.
In more detail, given two subsets $B=\{b_1,\dots,b_k\}$ and
$C=\{c_1,\dots,c_k\}$ of $\{1,\dots,n\}$ we will consider
the corresponding column-minor
\ben
M^B_C=\sum_{\si\in\Sym_k}\sgn\si\cdot
M_{b_{\si(1)}c_1}\dots M_{b_{\si(k)}c_k}.
\een
By \cite[Proposition~18]{cfr:ap} (see also \cite[Proposition~2.1]{mr:mm})
we have
\ben
A^{(m)} M_1\dots M_m = A^{(m)} M_1\dots M_m\tss A^{(m)},
\een
which implies
\beql{vpkm}
\tr_{1,\dots,m}\ts A^{(m)} M_1\dots M_m=\sum_{I,\ts |I|=m} M^I_I,
\eeq
summed over the subsets $I=\{i_1,\dots,i_m\}$ with $i_1<\dots<i_m$.
By Theorem~\ref{thm:comsera}, the algebra $\Ac_{\mu}$
is generated by the coefficients $\phi^{(k)}_{m}$ of the
constant term of the differential operator,
\ben
\phi^{\tss(0)}_{m}z^{-m}
+\dots+\phi^{\tss(m-1)}_{m}z^{-1}
+\phi^{\tss(m)}_{m}=\sum_{I,\ts |I|=m} M^I_I\ts 1,
\een
assuming that $\di_z\ts 1=0$. This implies the formula
\beql{defvp}
\phi^{(k)}_{m}=z^{m-k}\ts\sum_{I,\ts |I|=m}
\sum_{\atopn{B,C\subset I}{|B|=|C|=k}}\ts
\sgn\si\cdot{\mu\tss}^B_C\ts
\big[{-}\di_z+Ez^{-1}\big]^{I\setminus B}_{I\setminus C}\ts 1,
\eeq
where $\si$ denotes the permutation of the set $I$ given by
\ben
\si=\binom{B,\ I\setminus B}{C,\ I\setminus C}
=\binom{b_1,\dots,b_k,i_1,\dots,\wh b_1,\dots,\wh b_k,\dots,i_m}
{c_1,\dots,c_k,i_1,\dots,\wh c_1,\dots,\wh c_k,\dots,i_m},
\een
and we assume that $b_1<\dots<b_k$ and $c_1<\dots<c_k$ for the respective
elements of the subsets $B$ and $C$ in $I$.

For each $l=1,\dots,n$ introduce the polynomial in a variable $t$
with coefficients in $\Ac_{\mu}$ by
\beql{polphi}
\Phi_l(t,\mu)=\phi^{(0)}_{l}(\mu)\tss t^{\ts n-l}+\phi^{(1)}_{l+1}(\mu)\tss t^{\ts n-l-1}
+\dots+\phi^{(n-l)}_{n}(\mu),
\eeq
where the elements $\phi^{(k)}_{m}=\phi^{(k)}_{m}(\mu)$ are defined in \eqref{defvp}
and we indicated dependence of $\mu$.
The coefficients of $\Phi_l(t,\mu)$ are the elements of the $l$-th row
of the diagram $\Ga$; see \eqref{Ga}.

\ble\label{lem:relshi}
For any $a\in\CC$ we have the relation
\ben
\Phi_l(t,\mu+a\tss 1)=\Phi_l(t+a,\mu).
\een
\ele

\bpf
We have
\ben
\tr_{1,\dots,m}\ts A^{(m)} (a+M_1)\dots
(a+M_m)=\sum_{p=0}^m a^p \sum_{i_1<\dots<i_{m-p}}
\tr_{1,\dots,m}\ts A^{(m)} M_{i_1}\dots
M_{i_{m-p}}.
\een
Furthermore,
$A^{(m)}=\sgn p\cdot A^{(m)}\tss P$ for any $p\in\Sym_m$, where $P$ denotes
the image of $p$ in the algebra \eqref{tenprka}
under the action of $\Sym_m$. Hence, applying conjugations by appropriate elements
$P$ and using the cyclic property of trace, we can write the expression as
\ben
\sum_{p=0}^m \binom{m}{p}\tss a^p\ts
\tr_{1,\dots,m}\ts A^{(m)} M_1\dots M_{m-p}.
\een
The partial trace
of the anti-symmetrizer over the $m$-th copy of $\End\CC^n$ is found by
\beql{partr}
\tr_{m}\ts
A^{(m)}=\frac{n-m+1}{m} \ts A^{(m-1)}
\eeq
which implies
\ben
\tr_{m-p+1,\dots,m}\ts A^{(m)}=\frac{(n-m+p)!\ts(m-p)!}{(n-m)!\ts m!}\ts A^{(m-p)}.
\een
Hence,
\ben
\tr_{1,\dots,m}\ts A^{(m)} (a+M_1)\dots
(a+M_m)=\sum_{p=0}^m \binom{n-m+p}{p}\tss a^p\ts
\tr_{1,\dots,m-p}\ts A^{(m-p)} M_1\dots
M_{m-p}.
\een
Now equate the constant terms of the differential operators
on both sides and take the
coefficients of $z^{-m+k}$ to get the relation
\ben
\phi^{(k)}_{m}(\mu+a\tss 1)=\sum_{p=0}^k \binom{n-m+p}{p}\tss a^p\ts
\phi^{(k-p)}_{m-p}(\mu).
\een
Therefore, for the polynomial $\Phi_l(t,\mu+a\tss 1)$ we find
\ben
\bal
\Phi_l(t,\mu+a\tss 1)&=\sum_{k=0}^{n-l} \phi^{(k)}_{l+k}(\mu+a\tss 1)\tss t^{\ts n-l-k}
=\sum_{k=0}^{n-l} \ts\sum_{p=0}^k \binom{n-l-k+p}{p}\tss a^p\ts
\phi^{(k-p)}_{l+k-p}(\mu)\tss t^{\ts n-l-k}\\
{}&=\sum_{p=0}^{n-l}a^p\ts \sum_{r=0}^{n-l-p} \binom{n-l-r}{p}\ts
\phi^{(r)}_{l+r}(\mu)\tss t^{\ts n-l-p-r},
\eal
\een
which coincides with
\ben
\sum_{p=0}^{n-l}\frac{a^p}{p!}\ts \Big(\frac{d}{d\tss t}\Big)^p\ts
\Phi_l(t,\mu)=\Phi_l(t+a,\mu),
\een
as claimed.
\epf

\ble\label{lem:zecor}
Suppose that $\mu$ has the form of a block-diagonal matrix
\beql{mucan}
\mu=\begin{bmatrix}J_{\al}(0)&\text{\rm O}\ts\\
\text{\rm O}&\wt\mu\ts
\end{bmatrix},
\eeq
where $J_{\al}(0)$ is the nilpotent Jordan matrix associated
with a diagram $\al=(\al_1,\al_2,\dots)$ and $\wt\mu$ is an arbitrary square matrix
of size $q$ such that $|\al|+q=n$. Then for any $l\geqslant 1$ we have
\ben
\phi^{(k)}_{l+k}=0\qquad \text{for all}\quad
n-l-\de_l+1\leqslant k\leqslant n-l,
\een
where $\de_l=\al_{l+1}+\al_{l+2}+\dots$ is the number of boxes of $\al$
below its row $l$.
\ele

\bpf
The generator $\phi^{(k)}_{l+k}$ is found by \eqref{defvp}
for $m=l+k$.
The internal sum is a linear combination of
$k\times k$ minors of the matrix $\mu$ satisfying the condition that
the union $B\cup C$ of the row and column indices of each minor is a set of size
not exceeding $k+l$. On the other hand, with the given condition on $k$,
the minor ${\mu\tss}^B_C$ can be nonzero only if
the union of row and column indices is of the size at least $k+l+1$. Indeed,
this follows from the observation that if $p$ is a positive integer,
then any nonzero $p\times p$ minor of a nilpotent Jordan block
has the property that the minimal possible
size of the union of its row and column
indices is $p+1$. However,
the condition
$k\geqslant n-l-\de_l+1$ means that
$k\geqslant \al_1+\dots+\al_l-l+1+q$. Therefore,
a nonzero $k\times k$ minor must involve at least $l+1$ Jordan blocks.
\epf

In the following we use the notation of the Main Theorem. In addition, for
each diagram $\al^{(i)}$ we denote by $\de^{(i)}_l$ the corresponding
parameter $\de_l$, as defined in Lemma~\ref{lem:zecor}, so that
for the number $\ga_l$ defined in \eqref{gal} we have
\ben
\ga_l=\sum_{i=1}^r \de^{(i)}_l.
\een

\bco\label{cor:vanish}
The polynomial $\Phi_l(t,\mu)$ admits the factorization
\ben
\Phi_l(t,\mu)=(t+\la_1)^{\de^{(1)}_l}\dots (t+\la_r)^{\de^{(r)}_l}\ts \wt\Phi_l(t,\mu)
\een
for a certain polynomial $\wt\Phi_l(t,\mu)$ in $t$.
\eco

\bpf
The algebra $\Ac_{\mu}$ is known to depend only on
the adjoint orbit of $\mu$; see \cite{fft:gm}. More precisely, as we can see
from formulas \eqref{vpkm}, the elements $\phi^{(k)}_m$ are unchanged under
the simultaneous replacements $\mu\mapsto g\tss\mu\tss g^{-1}$
and $E\mapsto g\tss E\tss g^{-1}$ for $g\in \GL_n$. This implies that
$\Ac_{g\tss\mu\tss g^{-1}}$ can be identified with the algebra $\Ac_{\mu}$
associated with the image of $\U(\gl_n)$ under the automorphism
sending $E$ to $g\tss E\tss g^{-1}$.

For any $i\in\{1,\dots,r\}$
the Jordan canonical form of $\mu-\la_i 1$ is a matrix of the form
\eqref{mucan}, where $\al=\al^{(i)}$.
By Lemma~\ref{lem:zecor}, the polynomial $\Phi_l(t,\mu-\la_i 1)$
is divisible by $t^{\de^{(i)}_l}$. Hence, by
Lemma~\ref{lem:relshi}, the polynomial
$
\Phi_l(t,\mu)=\Phi_l(t+\la_i,\mu-\la_i 1)
$
is divisible by $(t+\la_i)^{\de^{(i)}_l}$.
\epf

We can now complete the proof of the Main Theorem. First, Corollary~\ref{cor:vanish}
implies that for any $l=1,\dots,n$ the generators $\phi^{(k)}_{l+k}$
with $n-l-\ga_l+1\leqslant k\leqslant n-l$ are linear combinations
of those generators with $k=0,1,\dots,n-l-\ga_l$. Therefore,
the elements $\phi^{(k)}_{l+k}$
corresponding to the boxes of the skew diagram $\Ga/\ga$
generate the algebra $\Ac_{\mu}$.
It remains to verify that these generators
are algebraically independent.

Consider the elements $\overline\phi^{\ts(k)}_{m}\in\Sr(\gl_n)$
which are defined by
\beql{bardefvp}
\overline\phi^{\ts(k)}_{m}=\sum_{I,\ts |I|=m}
\sum_{\atopn{B,C\subset I}{|B|=|C|=k}}\ts
\sgn\si\cdot{\mu\tss}^B_C\ts {E\tss}^{I\setminus B}_{I\setminus C},
\eeq
with the notation as in \eqref{defvp}, where the
entries of the matrix $E$ are now regarded as elements of the symmetric algebra
$\Sr(\gl_n)$. Equivalently, the elements $\overline\phi^{\ts(k)}_{m}$ are found by
\beql{phiz}
\tr_{1,\dots,m}\ts A^{(m)} \big(\mu^{}_1+E^{}_1z^{-1}\big)\dots
\big(\mu^{}_m+E^{}_mz^{-1}\big)
=\overline\phi^{\ts(0)}_{m}z^{-m}
+\dots+\overline\phi^{\ts(m-1)}_{m}z^{-1}
+\overline\phi^{\ts(m)}_{m}.
\eeq
They are generators of the subalgebra $\overline\Ac_{\mu}$.
The arguments of this section (including Lemmas~\ref{lem:relshi}, \ref{lem:zecor}
and Corollary~\ref{cor:vanish})
applied to these generators instead of the
$\phi^{(k)}_{m}$ show that the elements $\overline\phi^{\ts(k)}_{m}$
corresponding to the boxes of the skew diagram $\Ga/\ga$
generate the algebra $\overline\Ac_{\mu}$. Furthermore, we have
the following.

\ble\label{lem:algin}
The generators $\overline\phi^{\ts(k)}_{m}$ of the subalgebra $\overline\Ac_{\mu}$
corresponding to the boxes of the skew diagram $\Ga/\ga$
are algebraically independent.
\ele

\bpf
Regarding the elements $\overline\phi^{\ts(k)}_{m}$ as polynomials in the
variables $E_{ij}$, we will see that their differentials $d\ts\overline\phi^{\ts(k)}_{m}$
are linearly independent at a certain point.
Since these elements generate $\overline\Ac_{\mu}$, the linear span
of the differentials $d\ts\overline\phi^{\ts(k)}_{m}$
at any point coincides with the linear span
of all differentials
\ben
d\ts\overline\Ac_{\mu}=\text{span of}\ \{d\phi\ |\ \phi\in \overline\Ac_{\mu}\}.
\een
On the other hand, Bolsinov's criterion \cite[Theorem~3.2]{b:cf}
implies that the relation
\ben
\dim d\ts\overline\Ac_{\mu}=\text{\rm rank\ts}\gl_n+\frac12
\big(\hspace{-1pt}\dim\gl_n-\dim\gl_n^{\ts\mu}\big)
\een
holds at a certain regular point
if and only if the equality \eqref{beconj} holds for $\g=\gl_n$;
see also \cite[Theorem~2.7]{cm:ic}
for a concise exposition of this result. This equality does hold \cite{y:ic},
and so, to show that
the differentials $d\ts\overline\phi^{\ts(k)}_{m}$ of the generators
are linearly independent at a certain point,
we only
need to verify that the number of boxes of the skew diagram $\Ga/\ga$
coincides with
\ben
\text{\rm rank\ts}\gl_n+\frac12\big(\hspace{-1pt}\dim\gl_n-\dim\gl_n^{\ts\mu}\big)=
n+\frac12\big(n^2-\dim\gl_n^{\ts\mu}\big).
\een
Since $|\Ga|=n(n+1)/2$, the desired formula
is equivalent to the relation
\beql{relcen}
\dim\gl_n^{\ts\mu}=2\tss |\ga|+n.
\eeq
For the dimension of the centralizer we have
\ben
\dim\gl_n^{\ts\mu}=\sum_{i=1}^r\dim\gl_{n_i}^{\ts\mu^{(i)}},
\een
where $\mu^{(i)}$ denotes the direct sum of all Jordan blocks
of $\mu$ with the eigenvalue $\la_i$, and $n_i$ is the size of $\mu^{(i)}$.
Hence, by the definition of $\ga$, the verification of \eqref{relcen}
reduces to the case where $\mu$ has only one eigenvalue.
Let $\al_1\geqslant\dots\geqslant\al_s$ be the respective sizes
of the Jordan blocks of such matrix $\mu$. Then
$
\dim\gl_n^{\ts\mu}=\al_1+3\tss \al_2+\dots+(2\tss s-1)\tss\al_s,
$
while
\ben
|\ga|=\al_2+2\tss\al_3+\dots+(s-1)\tss \al_s\Fand n=\al_1+\dots+\al_s,
\een
thus implying \eqref{relcen}.
\epf

Now consider the generators $\phi^{(k)}_{m}$ of the algebra $\Ac_{\mu}$
associated with the boxes of the diagram $\Ga/\ga$.
By Lemma~\ref{lem:algin}, the corresponding
elements $\overline\phi^{\ts(k)}_{m}$ are nonzero, so that
the image of $\phi^{(k)}_{m}$
in the $(m-k)$-th component of $\gr\ts\U(\gl_n)\cong\Sr(\gl_n)$ coincides with
$\overline\phi^{\ts(k)}_{m}$. Moreover, the generators $\phi^{(k)}_{m}$
corresponding to the boxes of the diagram $\Ga/\ga$ are algebraically
independent. This completes the proof of the first part
of the Main Theorem, and the second part also follows.

\medskip

Finally, we will extend the first part of the Main Theorem by providing
some other families of algebraically independent generators of
the algebra $\Ac_{\mu}$.  To this end, introduce
the families $\overline\psi^{\ts(k)}_{m}$ and $\overline\theta^{\ts(k)}_{m}$
of generators of the algebra $\overline\Ac_{\mu}$
by the respective expansions
\beql{psiz}
\tr_{1,\dots,m}\ts H^{(m)} \big(\mu^{}_1+E^{}_1z^{-1}\big)\dots
\big(\mu^{}_m+E^{}_mz^{-1}\big)
=\overline\psi^{\ts(0)}_{m}z^{-m}
+\dots+\overline\psi^{\ts(m-1)}_{m}z^{-1}
+\overline\psi^{\ts(m)}_{m}
\eeq
and
\beql{thez}
\tr\ts\big(\mu+E\tss z^{-1}\big)^m=\overline\theta^{\ts(0)}_{m}z^{-m}
+\dots+\overline\theta^{\ts(m-1)}_{m}z^{-1}
+\overline\theta^{\ts(m)}_{m},
\eeq
where
\ben
E=\sum_{i,j=1}^n e_{ij}\ot E_{ij}\in \End\CC^n\ot \Sr(\gl_n),
\een
and extend the notation \eqref{matnota} to matrices $E$ and $\mu$.
The polynomials $\overline\phi_{m}(z)$, $\overline\psi_{m}(z)$ and
$\overline\theta_{m}(z)$ in $z^{-1}$ given by the respective expressions in
\eqref{phiz}, \eqref{psiz} and \eqref{thez} are related
by the classical MacMahon Master Theorem and Newton's identities:
\beql{mm}
\sum_{l=0}^m(-1)^l\ts \overline\phi_{l}(z)\ts \overline\psi_{m-l}(z)=0
\eeq
and
\beql{newton}
m\ts \overline\phi_{m}(z)=\sum_{l=1}^m (-1)^{l-1}\ts \overline\theta_{l}(z)
\ts \overline\phi_{m-l}(z)
\eeq
for $m\geqslant 1$.
Writing the relations \eqref{mm} and \eqref{newton} in terms of the coefficients
of the polynomials, we find that each of the generators
$\overline\psi^{\ts(k)}_{m}$ and $\overline\theta^{\ts(k)}_{m}$
with $m=1,\dots,n$ and $k=0,1,\dots,m-1$ will be presented in the form
\ben
c\cdot \overline\phi^{\ts(k)}_{m}+\text{linear combination of}\quad
\overline\phi^{\ts(k_1)}_{m_1}\dots \overline\phi^{\ts(k_s)}_{m_s},\quad s\geqslant 2,
\een
for a nonzero constant $c$,
where $m_1+\dots+m_s=m$ and $k_1+\dots+k_s=k$. As we pointed out above,
the elements $\overline\phi^{\ts(k)}_{m}$ corresponding to a certain row
of the diagram $\Ga$ are linear combinations of the elements of this row
in the skew diagram $\Ga/\ga$. This implies that each of the families
of generators $\overline\psi^{\ts(k)}_{m}$ and $\overline\theta^{\ts(k)}_{m}$
associated with the boxes of $\Ga/\ga$ as in \eqref{Ga}, is algebraically independent.
This leads to the following corollary,
where, as before, $\mu\in\gl_n$ is an arbitrary matrix.

\bco\label{cor:othgen}
The elements of each of the two families $\psi^{(k)}_{m\tss m}$ and $\theta^{\tss(k)}_{m\tss m}$
associated with the boxes
of the skew diagram $\Ga/\ga$ as in \eqref{Ga} are algebraically independent generators
of the algebra $\Ac_{\mu}$.
\qed
\eco

To construct two more families of generators of the algebra $\Ac_{\mu}$,
define the elements $\vp^{(k)}_m, \psi^{(k)}_m\in\U(\gl_n)$ by the expansions
\ben
\bal
\tr_{1,\dots,m}\ts A^{(m)} \big(\mu^{}_1+E^{}_1z^{-1}\big)\dots
\big(\mu^{}_m+E^{}_mz^{-1}\big)
{}&=\vp^{\tss(0)}_{m}z^{-m}
+\dots+\vp^{\tss(m-1)}_{m}z^{-1}
+\vp^{\tss(m)}_{m},\\[0.5em]
\tr_{1,\dots,m}\ts H^{(m)} \big(\mu^{}_1+E^{}_1z^{-1}\big)\dots
\big(\mu^{}_m+E^{}_mz^{-1}\big){}&=\psi^{\tss(0)}_{m}z^{-m}
+\dots+\psi^{\tss(m-1)}_{m}z^{-1}
+\psi^{\tss(m)}_{m}.
\eal
\een
It is easy to verify that each of the families
$\vp^{(k)}_m$ and $\psi^{(k)}_m$ with $m=1,\dots,n$ and $k=0,\dots,m-1$
generates the algebra $\Ac_{\mu}$. Indeed,
by Theorem~\ref{thm:comsera}, the algebra $\Ac_{\mu}$
is generated by the coefficients $\phi^{(k)}_{m}$ of the
constant term of the differential operator,
\begin{multline}
\tr_{1,\dots,m}\ts A^{(m)} \big({-}\di_z+\mu^{}_1+E^{}_1z^{-1}\big)\dots
\big({-}\di_z+\mu^{}_m+E^{}_mz^{-1}\big)\ts 1\\[0.7em]
{}=\phi^{\tss(0)}_{m}z^{-m}
+\dots+\phi^{\tss(m-1)}_{m}z^{-1}
+\phi^{\tss(m)}_{m}.
\non
\end{multline}
Hence,
$\phi^{(k)}_{m}$ is found as the coefficient of $z^{-m+k}$ in the expression
\ben
\sum_{i_1<\dots<i_k}
\sum_{j_1<\dots<j_{m-k}}
\ts\tr_{1,\dots,m}\ts A^{(m)} \mu^{}_{i_1}\dots \mu^{}_{i_k}
\big({-}\di_z+E^{}_{j_1}z^{-1}\big)\dots
\big({-}\di_z+E^{}_{j_{m-k}}z^{-1}\big)\ts 1,
\een
summed over disjoint subsets of indices $\{i_1,\dots,i_k\}$ and
$\{j_1,\dots,j_{m-k}\}$ of $\{1,\dots,m\}$. Therefore,
\ben
\phi^{(k)}_{m}=z^{m-k}\ts\binom{m}{k}\ts\tr_{1,\dots,m}
\ts A^{(m)} \mu^{}_1\dots \mu^{}_k
\big({-}\di_z+E^{}_{k+1}z^{-1}\big)\dots
\big({-}\di_z+E^{}_mz^{-1}\big)\ts 1.
\een
By calculating the partial trace
of the anti-symmetrizer with the use of \eqref{partr}, we get
\ben
\phi^{(k)}_{m}=\binom{m}{k}\ts\tr_{1,\dots,m}\ts A^{(m)} \mu^{}_1\dots \mu^{}_k
\ts E^{}_{k+1}\dots E^{}_m+\sum_{r=k+1}^{m-1} c_r\ts\tr_{1,\dots,r}\ts
A^{(r)} \mu^{}_1\dots \mu^{}_k
\ts E^{}_{k+1}\dots E^{}_r
\een
for certain constants $c_r$. The same argument applied to the expansion
defining the elements $\vp^{(k)}_{m}$
gives
\ben
\vp^{(k)}_{m}=\binom{m}{k}\ts\tr_{1,\dots,m}\ts A^{(m)} \mu^{}_1\dots \mu^{}_k
\ts E^{}_{k+1}\dots E^{}_m.
\een
This yields a triangular system of linear
relations
\ben
\phi^{(k)}_{m}=\vp^{(k)}_{m}+\sum_{r=k+1}^{m-1} c_r\ts \vp^{(k)}_{r}.
\een
Since $\phi^{(k)}_{k+1}=\vp^{(k)}_{k+1}$,
we can conclude that the elements $\vp^{(k)}_{m}$ are generators of $\Ac_{\mu}$.
The argument for the elements $\psi^{(k)}_{m}$ is quite similar.
Taking into account the properties of the elements
$\overline\phi^{\ts(k)}_{m}$ and $\overline\psi^{\ts(k)}_{m}$,
we come to another corollary.

\bco\label{cor:othgentwo}
The elements of each of the two families $\vp^{(k)}_{m}$ and $\psi^{(k)}_{m}$
associated with the boxes
of the skew diagram $\Ga/\ga$ as in \eqref{Ga} are algebraically independent generators
of the algebra $\Ac_{\mu}$.
\qed
\eco

\end{document}